\documentclass[12pt]{article}
\usepackage{amssymb}
\usepackage{amsmath}

\begin{document}

\begin{center}
\textbf{A generalization of Matkowski's fixed point theorem and Istr\u{a}%
\c{t}escu's fixed point theorem concerning convex contractions}

\bigskip

Radu MICULESCU and Alexandru MIHAIL

\bigskip
\end{center}

\textbf{Abstract}. {\small In this paper we obtain a generalization of
Matkowski's fixed point theorem and Istr\u{a}\c{t}escu's fixed point theorem
concerning convex contractions. More precisely, given a complete }$b${\small %
-metric space }$(X,d)${\small , we prove that every continuous function }$%
f:X\rightarrow X${\small \ is a Picard operator, provided that there exist }$%
m\in \mathbb{N}^{\ast }${\small \ and\ a comparison function }$\varphi $%
{\small \ such that }$d(f^{[m]}(x),f^{[m]}(y))\leq \varphi (\max
\{d(x,y),d(f(x),f(y)),...,d(f^{[m-1]}(x),f^{[m-1]}(y))\})${\small \ for all }%
$x,y\in X$.{\small \ In addition, we point out that if }$m=1${\small , the
continuity condition on }$f$ {\small is not necessary and consequently,
taking into account that a metric space is a }$b${\small -metric space, we
obtain a generalization of Matkowski's fixed point theorem. Moreover, we
prove that Istr\u{a}\c{t}escu's fixed point theorem concerning convex
contractions is a particular case of our result for }$m=2${\small . By
providing appropriate examples we show that the above mentioned two
generalizations are effective.}

\bigskip



\bigskip

\textbf{1.} \textbf{Introduction}

\bigskip

One branch of generalizations of the celebrated Banach-Cacciopoli-Picard
contraction principle is based on the replacement of the contractivity
condition imposed on the function $f:X\rightarrow X$, where $(X,d)$ is a
complete metric metric space, by a weaker one described by the inequality $%
d(f(x),f(y))\leq \varphi (d(x,y))$ for all $x,y\in X$, where $\varphi $ has
certain properties (see [8], [9], [17] and [23]). The result obtained by J.
Matkowski in this direction can be stated as follows: \textit{Every }$%
\varphi $\textit{-contraction }$f:X\rightarrow X$\textit{, where }$(X,d)$%
\textit{\ is a complete metric space, is a Picard operator.}

The notion of $b$-metric space was introduced by I. A. Bakhtin [4] and S.
Czerwik (see [11] and [12]) in connection with some problems concerning the
convergence of measurable functions with respect to measure. Among the fixed
point results in the framework of $b$-metric spaces obtained in the last
period (see, for example, [2], [6], [7], [10], [18], [19], [22], [25], [28]
and the references therein, to mention just the most recent ones) we point
out a Matkowski type fixed point result (in the framework of such a space $X$
with the property that the $b$-metric is a continuous functional on $X\times
X$ and for a $b$-comparison function $\varphi $) which is due to M. P\u{a}%
curar (Berinde) (see [20] and [21]).

Motivated by the fact that in some situations there is no need to use the
entire force of metric requirements in the proof of certain fixed point
theorems, J. Jachymski, J. Matkowski and T. \'{S}wi\c{a}tkowski [16]
obtained a generalization of Matkowski's result for the class of semi-metric
spaces satisfying the so called JMS condition (see [3]). Note that this
class is larger that the one of $b$-metric spaces. Another proof of this
result (which is based on a monotone principle of fixed points -see [26]-)
was given by M. Taskovi\'{c} [27]. By defining the notion of $\tau $%
-distance function in a general topological space $(X,\tau )$, M. Aauri and
D. El. Montawokil [1] obtained a generalization of the result due to
Jachymski, Matkowski and \'{S}wi\c{a}tkowski. A related result (namely an
extension of the Matkowski fixed point theorem in the framework of complete
and regular semi-metric spaces) is given in [5].

In an attempt to study if there exist contraction-type conditions that do no
imply the contraction condition and for which the existence and uniqueness
of the fixed point are assured, V. Istr\u{a}\c{t}escu introduced and studied
the convex contraction condition (see [13], [14] and [15]). A continuous
function $f:(X,d)\rightarrow (X,d)$, where $(X,d)$ is a complete metric
space, is called convex contraction if there exist $a,b\in (0,1)$ such that $%
a+b<1$ and $d(f^{[2]}(x),f^{[2]}(y))\leq ad(f(x)),f(y))+bd(x,y)$ for all $%
x,y\in X$. Istr\u{a}\c{t}escu proved that any convex contraction has a
unique fixed point $\alpha \in X$ (and $\underset{n\rightarrow \infty }{\lim 
}f^{[n]}(x)=\alpha $ for every $x\in X$) and provided an example of convex
contraction which is not contraction.

In this paper we prove the following result: Every continuous function $%
f:X\rightarrow X$, where $(X,d)$\ is a complete $b$-metric space, is a
Picard operator, provided that there exist $m\in \mathbb{N}^{\ast }$ and\ a
comparison function $\varphi $ such that $d(f^{[m]}(x),f^{[m]}(y))\leq
\varphi (\max \{d(x,y),d(f(x),f(y)),...,d(f^{[m-1]}(x),f^{[m-1]}(y))\})$ for
all $x,y\in X$. Moreover, we show that if $m=1$, then the continuity of $f$
is not necessary and therefore we obtain Matkowski's fixed point theorem as
a particular case (by taking into account the fact that every metric space
is a $b$-metric space). In addition, we proved that Istr\u{a}\c{t}escu's
fixed point theorem concerning convex contractions is a particular case of
our result for $m=2$. We provide two examples to justify that the above
mentioned two generalizations are effective.

\newpage

\textbf{2.} \textbf{Preliminaries}

\bigskip

In this section we recall some basic facts that will be used in the sequel.

\bigskip

\textbf{Definition 2.1.} \textit{Given a nonempty set }$X$\textit{\ and a
real number }$s\in \lbrack 1,\infty )$\textit{, a function }$d:X\times
X\rightarrow \lbrack 0,\infty )$\textit{\ is called }$b$\textit{-metric if
it satisfies the following properties:}

\textit{i) }$d(x,y)=0$\textit{\ if and only if }$x=y$\textit{;}

\textit{ii) }$d(x,y)=d(y,x)$\textit{\ for all }$x,y\in X$\textit{;}

\textit{iii) }$d(x,y)\leq s(d(x,z)+d((z,y))$\textit{\ for all }$x,y,z\in X$%
\textit{.}

\textit{The pair }$(X,d)$\textit{\ is called }$b$\textit{-metric space (with
constant }$s$\textit{).}

\bigskip

Besides the classical spaces $l^{p}(\mathbb{R})$ and $L^{p}[0,1]$, where $%
p\in (0,1)$, more examples of $b$-metric spaces could be found in [6], [11]
and [12].

\bigskip

\textbf{Remark 2.1}. \textit{Every metric space is a }$b$\textit{-metric
space (with constant }$1$\textit{). There exist }$b$\textit{-metric spaces
which are not metric spaces (see, for example, }[10] \textit{or} [19]\textit{%
).}

\bigskip

\textbf{Definition 2.2.} \textit{A sequence} $(x_{n})_{n\in \mathbb{N}}$ 
\textit{of elements from a }$b$\textit{-metric space }$(X,d)$\textit{\ is
called:}

\textit{- convergent if there exists }$l\in \mathbb{R}$\textit{\ such that }$%
\underset{n\rightarrow \infty }{\lim }d(x_{n},l)=0$;

\textit{- Cauchy if }$\underset{m,n\rightarrow \infty }{\lim }%
d(x_{m},x_{n})=0$\textit{, i.e. for every }$\varepsilon >0$\textit{\ there
exists} $n_{\varepsilon }\in \mathbb{N}$ \textit{such that }$%
d(x_{m},x_{n})<\varepsilon $\textit{\ for all} $m,n\in \mathbb{N}$, $m,n\geq
n_{\varepsilon }$.

\textit{The} $b$\textit{-metric space }$(X,d)$ \textit{is called complete if
every Cauchy sequence of elements from }$(X,d)$\textit{\ is convergent.}

\bigskip

\textbf{Remark 2.2.} \textit{A }$b$\textit{-metric space can be endowed with
the topology induced by its convergence.}

\bigskip

For a function $f:X\rightarrow X$ and $n\in \mathbb{N}$, by $f^{[n]}$ we
mean the composition of $f$ by itself $n$ times.

\bigskip

\textbf{Definition 2.3}. \textit{A function} $\varphi :[0,\infty
)\rightarrow \lbrack 0,\infty )$ \textit{is called a comparison function if:}

\textit{i) }$\varphi $\textit{\ is increasing;}

\textit{ii) }$\underset{n\rightarrow \infty }{\lim }\varphi ^{\lbrack
n]}(r)=0$\textit{\ for every }$r\in \lbrack 0,\infty )$\textit{.}

\bigskip

\textbf{Remark 2.3}. \textit{Every comparison function }$\varphi $\textit{\
has the property that}\linebreak \textit{\ }$\varphi (0)=0$ \textit{and} $%
\varphi (r)<r$\textit{\ for every }$r\in (0,\infty )$\textit{.}

\bigskip

\textbf{Definition 2.4.} \textit{Given a }$b$\textit{-metric space }$(X,d)$%
\textit{\ and a comparison function }$\varphi $\textit{, a function }$%
f:X\rightarrow X$\textit{\ is called }$\varphi $\textit{-contraction if }$%
d(f(x),f(y))\leq \varphi (d(x,y))$\textit{\ for all }$x,y\in X$\textit{.}

\bigskip

\textbf{Definition 2.5. }\textit{Given a }$b$\textit{-metric space }$(X,d)$, 
\textit{a function} $f:X\rightarrow X$\textit{\ is called Picard operator if
there exists a unique fixed point }$\alpha $ \textit{of }$f$\textit{\ and
the sequence }$(f^{[n]}(x))_{n\in \mathbb{N}}$\textit{\ is convergent to }$%
\alpha $ \textit{for every }$x\in X$.

\bigskip

\textbf{Lemma 2.1}. \textit{Given\ a }$b$\textit{-metric space }$(X,d)$%
\textit{\ with constant }$s\geq 1$\textit{, the }\linebreak \textit{%
inequality }$d(x_{0},x_{p})\leq \underset{i=1}{\overset{p}{\sum }}%
s^{i}d(x_{i-1},x_{i})$\textit{\ is valid for every }$p\in \mathbb{N}$\textit{%
\ and for all }$x_{0},x_{1},...,x_{p}\in X$.

\textit{Proof}. We have%
\begin{equation*}
d(x_{0},x_{p})\leq sd(x_{0},x_{1})+sd(x_{1},x_{p})\leq
sd(x_{0},x_{1})+s^{2}d(x_{1},x_{2})+s^{2}d(x_{2},x_{p})\leq
\end{equation*}%
\begin{equation*}
\leq
sd(x_{0},x_{1})+s^{2}d(x_{1},x_{2})+s^{3}d(x_{2},x_{3})+...+s^{p-1}d(x_{p-2},x_{p-1})+s^{p-1}d(x_{p-1},x_{p})\leq
\end{equation*}%
\begin{equation*}
\leq
sd(x_{0},x_{1})+s^{2}d(x_{1},x_{2})+...+s^{p-1}d(x_{p-2},x_{p-1})+s^{p}d(x_{p-1},x_{p})=%
\underset{i=1}{\overset{p}{\sum }}s^{i}d(x_{i-1},x_{i})\text{,}
\end{equation*}%
for every\textit{\ }$p\in \mathbb{N}$\ and all\textit{\ }$%
x_{0},x_{1},...,x_{p}\in X$. $\square $

\bigskip

\textbf{3.} \textbf{The main result}

\bigskip

\textbf{Theorem 3.1}. \textit{Every continuous function }$f:X\rightarrow X$%
\textit{, where }$(X,d)$\textit{\ is a complete }$b$\textit{-metric space,
is a Picard operator, provided that there exist }$m\in \mathbb{N}^{\ast }$ 
\textit{and\ a comparison function }$\varphi $ \textit{such that the
following inequality: }%
\begin{equation*}
d(f^{[m]}(x),f^{[m]}(y))\leq \varphi (\max
\{d(x,y),d(f(x),f(y)),...,d(f^{[m-1]}(x),f^{[m-1]}(y))\})\text{,}
\end{equation*}%
\textit{is valid for all} $x,y\in X$.

\textit{Proof}. Let us denote the constant of the $b$-metric space $(X,d)$
by $s$.

In the sequel, for $x,y\in X$ and $n\in \mathbb{N}$, we adopt the following
notations:

i) $x_{n}:=f^{[n]}(x)$ and $y_{n}:=f^{[n]}(y)$;

ii) $M_{n}(x,y):=\max
\{d(x_{n},y_{n}),d(x_{n+1},y_{n+1}),..,d(x_{n+m-1},y_{n+m-1})\}$.

\medskip

\textbf{Claim 1}. \textit{The sequence} $(M_{n}(x,y))_{n\in \mathbb{N}}$ 
\textit{is decreasing.}

\textit{Justification of Claim 1}. Taking into the inequality from
hypothesis\linebreak\ $x=x_{n}$ and $y=y_{n}$ we get that $%
d(x_{m+n},y_{m+n})\leq \varphi (M_{n}(x,y))\overset{\text{Remark 2.3}}{\leq }%
M_{n}(x,y)$, so $M_{n+1}(x,y)\leq M_{n}(x,y)$ for every $n\in \mathbb{N}$.

\medskip

\textbf{Claim 2}. $\underset{n\rightarrow \infty }{\lim }d(x_{n},y_{n})=0$.

\textit{Justification of Claim 2}. For every $i\in \{0,1,...,m-1\}$ and $%
n\in \mathbb{N}$, taking $x=x_{n+i}$ and $y=y_{n+i}$ into the inequality
from hypothesis, we get that $d(x_{m+n+i},y_{m+n+i})\leq \varphi
(M_{n+i}(x,y))\overset{\text{Claim 1}}{\leq }\varphi (M_{n}(x,y))$, so $%
M_{m+n}(x,y)\leq \varphi (M_{n}(x,y))$. Using the mathematical induction
method, we obtain that $M_{n+km}(x,y)\leq \varphi ^{\lbrack k]}(M_{n}(x,y))$
for all $k,n\in \mathbb{N}$ and since $\underset{k\rightarrow \infty }{\lim }%
\varphi ^{\lbrack k]}(M_{n}(x,y))=0$ we infer that $\underset{k\rightarrow
\infty }{\lim }M_{n+km}(x,y)=0$. Using Claim 1 we deduce that\linebreak\ $%
\underset{n\rightarrow \infty }{\lim }M_{n}(x,y)=0$ and since $%
d(x_{n},y_{n})\leq M_{n}(x,y)$ for every $n\in \mathbb{N}$, we conclude that 
$\underset{n\rightarrow \infty }{\lim }d(x_{n},y_{n})=0$.

\medskip

By taking $y=f(x)$, from Claim 2, we obtain that%
\begin{equation}
\underset{n\rightarrow \infty }{\lim }d(x_{n},x_{n+1})=0\text{.}  \tag{1}
\end{equation}

The above inequality assures us that there exits $n_{0}\in \mathbb{N}$ such
that%
\begin{equation}
d(x_{n},x_{n+1})\leq 1\text{,}  \tag{2}
\end{equation}%
for every $n\in \mathbb{N}$, $n\geq n_{0}$.

\textbf{Claim 3}. \textit{The sequence }$(x_{n})_{n\in \mathbb{N}}$\textit{\
is Cauchy.}

\textit{Justification of Claim 3}. Let us suppose, by reductio ad absurdum,
that $(x_{n})_{n\in \mathbb{N}}$\textit{\ }is not Cauchy. Then there exists $%
\varepsilon _{0}>0$ having the property that for every $k\in \mathbb{N}$
there exist $m_{k},n_{k}\in \mathbb{N}$, $m_{k},n_{k}>k$ such that $%
d(x_{n_{k}},x_{m_{k}})\geq \varepsilon _{0}$. Hence we get two subsequences $%
(x_{n_{k}})_{k\geq 1}$ and $(x_{m_{k}})_{k\geq 1}$ of $(x_{n})_{n\in \mathbb{%
N}}$ satisfying the following properties:

a) $n_{1}\geq n_{0}$;

b) $d(x_{n_{k}},x_{m_{k}})\geq \varepsilon _{0}$;

c) $n_{k}<m_{k}=\min \{n_{k}+p\mid p\in \mathbb{N}$ and $%
d(x_{n_{k}},x_{n_{k}+p})\geq \varepsilon _{0}\}$.

In the sequel we adopt the following notation: $C:=s^{3}(2\frac{s^{m}-1}{s-1}%
+\varepsilon _{0}+1)$.

Let us note that%
\begin{equation}
d(x_{n_{k}},x_{m_{k}})<s(\varepsilon _{0}+1)\leq C\text{,}  \tag{3}
\end{equation}%
for every $k\in \mathbb{N}$.

Indeed, we have $d(x_{n_{k}},x_{m_{k}})\leq
sd(x_{n_{k}},x_{m_{k}-1})+sd(x_{m_{k}-1},x_{m_{k}})\overset{c)\&(2)}{\leq }%
s(\varepsilon _{0}+1)$ for every $k\in \mathbb{N}$.

Moreover we have%
\begin{equation}
d(x_{n_{k}+i},x_{m_{k}+i})\leq C\text{,}  \tag{4}
\end{equation}%
for every $k\in \mathbb{N}$ and every $i\in \{1,2,...,m-1\}$.

Indeed, we have%
\begin{equation*}
d(x_{n_{k}+i},x_{m_{k}+i})\leq
s^{2}(d(x_{n_{k}+i},x_{n_{k}})+d(x_{n_{k}},x_{m_{k}})+d(x_{m_{k}},x_{m_{k}+i}))\leq
\end{equation*}%
\begin{equation*}
\overset{\text{Lemma 2.1}}{\leq }s^{2}(\underset{l=1}{\overset{i}{\sum }}%
s^{l}d(x_{n_{k}+l-1},x_{n_{k}+l})+d(x_{n_{k}},x_{m_{k}})+\underset{l=1}{%
\overset{i}{\sum }}s^{l}d(x_{m_{k}+l-1},x_{m_{k}+l}))\leq
\end{equation*}%
\begin{equation*}
\overset{(2)}{\leq }s^{2}(2\underset{l=1}{\overset{i}{\sum }}%
s^{l}+d(x_{n_{k}},x_{m_{k}}))\overset{(3)}{\leq }s^{2}(2s\frac{s^{i}-1}{s-1}%
+s(\varepsilon _{0}+1))\leq C\text{,}
\end{equation*}%
for every $k\in \mathbb{N}$ and every $i\in \{1,2,...,m-1\}$.

Inequalities $(3)$ and $(4)$ guarantee that%
\begin{equation}
M_{0}(x_{n_{k}},x_{m_{k}})\leq C\text{,}  \tag{5}
\end{equation}%
for every $k\in \mathbb{N}$.

As $\varphi $\textit{\ }is a comparison function, $\underset{n\rightarrow
\infty }{\lim }\varphi ^{\lbrack n]}(C)=0$, so there exists $p\in \mathbb{N}$
such that 
\begin{equation}
s^{2}\varphi ^{\lbrack p]}(C)<\frac{\varepsilon _{0}}{3}\text{.}  \tag{6}
\end{equation}

Now we choose $\varepsilon _{1}>0$ such that%
\begin{equation}
2\varepsilon _{1}s^{3}\frac{s^{pm}-1}{s-1}<\frac{\varepsilon _{0}}{3}\text{.}
\tag{7}
\end{equation}

Moreover, taking into account $(1)$, there exists $n_{\varepsilon _{1}}\in 
\mathbb{N}$ such that%
\begin{equation}
d(x_{n},x_{n+1})\leq \varepsilon _{1}\text{,}  \tag{8}
\end{equation}%
for every $n\in \mathbb{N}$, $n\geq n_{\varepsilon _{1}}$.

Then, for $k\in \mathbb{N}$ such that $n_{k}>n_{\varepsilon _{1}}$, we have:%
\begin{equation*}
\varepsilon _{0}\overset{b)}{\leq }d(x_{n_{k}},x_{m_{k}})\leq
s^{2}d(x_{n_{k}},x_{n_{k}+pm})+s^{2}d(x_{n_{k}+pm},x_{m_{k}+pm})+s^{2}d(x_{m_{k}+pm},x_{m_{k}})\leq
\end{equation*}%
\begin{equation*}
\overset{\text{Lemma 2.1}}{\leq }s^{2}\varphi ^{\lbrack
p]}(M_{0}(x_{n_{k}},x_{m_{k}}))+s^{2}\underset{i=1}{\overset{pm}{\sum }}%
s^{i}(d(x_{n_{k}+i-1},x_{n_{k}+i})+d(x_{m_{k}+i-1},x_{m_{k}+i}))\leq
\end{equation*}%
\begin{equation*}
\overset{(8)}{\leq }2\varepsilon _{1}s^{2}\underset{i=1}{\overset{pm}{\sum }}%
s^{i}+s^{2}\varphi ^{\lbrack p]}(M_{0}(x_{n_{k}},x_{m_{k}}))\overset{(5)}{%
\leq }2\varepsilon _{1}s^{3}\frac{s^{pm}-1}{s-1}+s^{2}\varphi ^{\lbrack
p]}(C)\overset{(6)\&(7)}{\leq }\frac{2\varepsilon _{0}}{3}\text{.}
\end{equation*}

This contradiction closes the justification of the claim.

\medskip

Since $(X,d)$ is a complete $b$-metric space, Claim 3 assures us that there
exists $\alpha \in X$ such that 
\begin{equation}
\underset{n\rightarrow \infty }{\lim }x_{n}=\alpha \text{.}  \tag{9}
\end{equation}

As $f$ is continuous, from $(9)$ we infer that $f(\alpha )=\underset{%
n\rightarrow \infty }{\lim }f(x_{n})=\underset{n\rightarrow \infty }{\lim }%
x_{n+1}$, so%
\begin{equation}
\underset{n\rightarrow \infty }{\lim }x_{n}=f(\alpha )\text{.}  \tag{10}
\end{equation}

Relations $(9)$ and $(10)$ imply that $f(\alpha )=\alpha $, i.e. $\alpha $
(which is the limit of the sequence $(f^{[n]}(x))_{n\in \mathbb{N}}$ with $x$
arbitrarily chosen in $X$) is a fixed point of $f$.

In addition, $\alpha $ is the unique fixed point of $f$.

Indeed, if $\beta \in X\smallsetminus \{\alpha \}$ would be a fixed point of 
$f$, then we arrive at the following contradiction: $0<d(\alpha ,\beta
)=d(f^{[m]}(\alpha ),f^{[m]}(\beta ))\leq $\linebreak $\varphi (d(\alpha
,\beta ),$ $d(f(\alpha ),f(\beta )),...,d(f^{[m]}(\alpha ),f^{[m]}(\beta
))=\varphi (d(\alpha ,\beta ))\overset{\text{Remark 2.3}}{<}d(\alpha ,\beta
) $. $\square $

\bigskip

\textbf{Corollary 3.1 (Matkowski's fixed point theorem in the framework of }$%
\mathbf{b}$\textbf{-metric spaces)}. \textit{Every }$\varphi $\textit{%
-contraction }$f:X\rightarrow X$\textit{, where }$(X,d)$\textit{\ is a
complete }$b$\textit{-metric space and }$\varphi $\textit{\ a comparison
function, is a Picard operator.}

\textit{Proof}. The proof is exactly as the one of the above result, except
the justification of $(10)$ which is the following one: because $%
d(x_{n+1},f(\alpha ))=d(f(x_{n}),f(\alpha ))\leq \varphi (d(x_{n},\alpha ))%
\overset{\text{Remark 2.3}}{\leq }d(x_{n},\alpha )$ for every $n\in \mathbb{N%
}$, using $(9)$, we come to the conclusion that $\underset{n\rightarrow
\infty }{\lim }x_{n}=f(\alpha )$. $\square $

\bigskip

\textbf{Example 3.1.} The function $f:[0,1]\rightarrow \lbrack 0,1]$ given
by $f(x)=\{%
\begin{array}{cc}
\frac{x}{4}\text{,} & x\in \lbrack 0,\frac{1}{2}) \\ 
\frac{x}{5}\text{,} & x\in \lbrack \frac{1}{2},1]%
\end{array}%
$ -see [13]- has the property that $d(f^{[2]}(x),f^{[2]}(y))\leq \varphi
(\max \{d(x,y),d(f(x),f(y))\})$ for all $x,y\in \lbrack 0,1]$, where $%
\varphi $ is the comparison function given by $\varphi (r)=\frac{1}{4}r$ for
every $r\in \lbrack 0,\infty )$, so it satisfies the hypothesis of Theorem
3.1. Since $f$ is not continuous, it does not satisfy the hypothesis of
Matkowski's Theorem. Thus \textit{Theorem 3.1 is an effective generalization
of Matkowski's fixed point theorem}.

\bigskip

\textbf{Corollary 3.2 (Istr\u{a}\c{t}escu's fixed point theorem concerning
convex contractions in the framework of }$\mathbf{b}$\textbf{-metric spaces)}%
. \textit{Given a complete }$b$\textit{-metric space }$(X,d)$, \textit{every
convex contraction }$f:X\rightarrow X$ \textit{(i.e. there exist }$a,b\in
(0,1)$\textit{\ such that }$a+b<1$\textit{\ and }$d(f^{[2]}(x),f^{[2]}(y))%
\leq ad(f(x)),f(y))+bd(x,y)$\textit{\ for all }$x,y\in X$\textit{) is a
Picard operator.}

\textit{Proof}. By considering the comparison function $\varphi :[0,\infty
)\rightarrow \lbrack 0,\infty )$ given by $\varphi (r)=(a+b)r$ for every $%
r\in \lbrack 0,\infty )$, we have 
\begin{equation*}
d(f^{[2]}(x),f^{[2]}(y))\leq ad(f(x)),f(y))+bd(x,y)\leq
\end{equation*}%
\begin{equation*}
\leq (a+b)\max \{d(x,y),d(f(x)),f(y))\}=\varphi (\max
\{d(x,y),d(f(x)),f(y))\})\text{,}
\end{equation*}%
for all $x,y\in X$. Now we just apply Theorem 3.1 for $m=2$. $\square $

\bigskip

\textbf{Example 3.2.} The function $f:[0,\frac{1}{2}]\rightarrow \lbrack 0,%
\frac{1}{2}]$, given by $f(x)=x-x^{2}$ for every $x\in \lbrack 0,\frac{1}{2}]
$, is not a convex contraction since if this is not the case, then there
exist $a,b\in (0,1)$ such that $a+b<1$ and $\left\vert
f^{[2]}(x)-f^{[2]}(y)\right\vert \leq a\left\vert f(x))-f(y)\right\vert
+b\left\vert x-y\right\vert $ for all $x,y\in \lbrack 0,\frac{1}{2}]$.
Picking $x\in \lbrack 0,\frac{1}{2}]$, one can easily check that the
sequence $(x_{n})_{n\in \mathbb{N}}$, given by $x_{n}=f^{[n]}(x)$ for every $%
n\in \mathbb{N}$, satisfies the following two properties: a) $x_{n}\leq
(a+b)^{[\frac{n}{2}]}$ for every $n\in \mathbb{N}$; b) $\underset{%
n\rightarrow \infty }{\lim }nx_{n}=1$. Consequently $nx_{n}\leq n(a+b)^{[%
\frac{n}{2}]}$ for every $n\in \mathbb{N}$ and by passing to limit as $n$
goes to $\infty $ we obtain the contradiction $1\leq 0$. Thus $f$ does not
satisfy the hypothesis of Istr\u{a}\c{t}escu's fixed point theorem
concerning convex contractions\textbf{.}

Considering the comparison function $\varphi :[0,\infty )\rightarrow \lbrack
0,\infty )$ given by $\varphi (x)=\{%
\begin{array}{cc}
x-x^{2}\text{,} & x\in \lbrack 0,\frac{1}{2}] \\ 
\frac{1}{4}\text{,} & x\in (\frac{1}{2},\infty )%
\end{array}%
$, one can easily check that $\left\vert f(x))-f(y)\right\vert \leq \varphi
(\left\vert x-y\right\vert )$ for all $x,y\in \lbrack 0,\frac{1}{2}]$.
Consequently $\left\vert f^{[2]}(x))-f^{[2]}(y)\right\vert \leq \varphi
(\left\vert f(x))-f(y)\right\vert )\leq \varphi (\max \{\left\vert
x-y\right\vert ,\left\vert f(x))-f(y)\right\vert \})$ for all $x,y\in
\lbrack 0,\frac{1}{2}]$. Hence $f$ satisfies the hypothesis of Theorem 3.1
for $m=2$.

Therefore \textit{Theorem 3.1 (for }$m=2$\textit{) is an effective
generalization of Istr\u{a}\c{t}escu's fixed point theorem concerning convex
contractions.}

\bigskip

\textbf{Remark 3.1}. \textit{Theorem 3.1 gives a partial answer, in the
framework of }$b$\textit{-metric spaces, to Problem 9.3.1 b) from }[24]%
\textit{.}

\bigskip

\textbf{References}

\bigskip

[1] M. Aauri and D. El. Montawokil, $\tau $-distance in a general
topological space $(X,\tau )$ with application to fixed point theory,
Southwest J. Pure Appl. Math., \textbf{2} (2003), 1-5.

[2] A. Aghajani, M. Abbas and J.R. Roshan, Common fixed point of generalized
weak contractive mappings in partially ordered $b$-metric spaces, Math.
Slovaca, \textbf{64} (2014), 941-960.

[3] S. Alshehri, I. Arandelovi\'{c} and N. Shahzad, Symmetric spaces and
fixed points of generalized contractions, Abstr. Appl. Anal., 2014, Article
ID 763547.

[4] I. A. Bakhtin, The contraction mapping principle in quasimetric spaces,
Funct. Anal., Unianowsk Gos. Ped. Inst., \textbf{30} (1989), 26-37.

[5] M. Bessenyei and Z. P\`{a}les, A contraction principle in semimetric
spaces, arXiv: 1401.1709.

[6] M. Bota, A. Moln\'{a}r and C. Varga, On Ekeland's variational principle
in $b$-metric spaces, Fixed Point Theory, \textbf{12} (2011), 21-28.

[7] M. Bota, V. Ilea, E. Karapinar, O. Mle\c{s}ni\c{t}e, On $\alpha _{\ast }$%
-$\varphi $-contractive multi-valued operators in $b$-metric spaces and
applications, Appl. Math. Inf. Sci., \textbf{9} (2015), 2611-2620.

[8] D.W. Boyd and J.S. Wong, On nonlinear contractions, Proc. Amer Math.
Soc., \textbf{20} (1969), 458-464.

[9] F.E. Browder, On the convergence of succesive approximations for
nonlinear functional equations, Nederl. Akad. Wet., Proc., Ser. A 71, 
\textbf{30} (1968), 27-35.

[10] C. Chifu and G. Petru\c{s}el, Fixed points for multivalued contractions
in $b$-metric spaces with applications to fractals, Taiwanese J. Math., 
\textbf{18} (2014), 1365-1375.

[11] S. Czerwik, Contraction mappings in $b$-metric spaces, Acta Math.
Inform. Univ. Ostraviensis, \textbf{1} (1993), 5-11.

[12] S. Czerwik, Nonlinear set-valued contraction mappings in $b$-metric
spaces, Atti Sem. Mat. Fis. Univ. Modena, \textbf{46} (1998), 263-276.

[13] V. Istr\u{a}\c{t}escu, Some fixed point theorems for convex contraction
mappings and convex nonexpansive mappings (I), Libertas Math., \textbf{1}
(1981), 151-164.

[14] V. Istr\u{a}\c{t}escu, Some fixed point theorems for convex contraction
mappings and mappings with convex diminishing diameters - I, Annali di Mat.
Pura Appl., \textbf{130} (1982), 89--104.

[15] V. Istr\u{a}\c{t}escu,\ Some fixed point theorems for convex
contraction mappings and mappings with convex diminishing diameters, II,
Annali di Mat. Pura Appl., \textbf{134} (1983), 327-362.

[16] J. Jachymski, J. Matkowski and T. \'{S}wi\c {a}tkowski, Nonlinear
contractions on semimetric spaces, J. Appl. Anal., \textbf{1} (1995),
125-134.

[17] J. Matkowski, Integrable solutions of functional equations,
Dissertations Math. (Rozprawy), 127 (1976).

[18] S. K. Mohanta, Some fixed point theorems using $wt$-distance in $b$%
-metric spaces, Fasc. Math., \textbf{54} (2015), 125-140.

[19] H.N. Nashine and Z. Kadelburg, Cyclic generalized $\varphi $%
-contractions in $b$-metric spaces and an application to integral equations,
Filomat, \textbf{28} (2014), 2047-2057.

[20] M. P\u{a}curar (Berinde), Iterative methods for fixed point
approximation, Ph. D. Thesis, Babe\c{s}-Bolyai University Cluj-Napoca,
Romania, 2009.

[21] M. P\u{a}curar (Berinde), A fixed point result for $\varphi $%
-contractions on $b$-metric spaces without the boundedness assumption, Fasc.
Math., \textbf{43} (2010), 127-137.

[22] J.R. Roshan, N. Hussain, S. Sedghi and N. Shobkolaei, Suzuki-type fixed
point results in $b$-metric spaces, Math. Sci. (Springer), \textbf{9}
(2015), 153--160.

[23] I. A. Rus, Generalized $\varphi $-contractions, Math., Rev. Anal. Num%
\'{e}r. Th\'{e}or. Approximation, Math., \textbf{47} (1982), 175-178.

[24] I. A. Rus, Generalized Contractions and Applications, Cluj University
Press, Cluj-Napoca, 2001.

[25] S. Shukla, Partial $b$-metric spaces and fixed point theorems,
Mediterr. J. Math., \textbf{11} (2014), 703-711.

[26] M. Taskovi\'{c}, A monotone principle of fixed points, Proc. Amer.
Math. Soc., \textbf{94} (1985), 427-432.

[27] M. Taskovi\'{c}, On a result of Jachymski, Matkowski, and \'{S}wi\c{a}%
tkowski, Math. Morav., \textbf{16} (2012), 33-35.

[28] H. Yingtaweesittikul, Suzuki type fixed point for generalized
multi-valued mappings in $b$-metric spaces, Fixed Point Theory Appl., 2013,
2013:215.

\bigskip

University of Bucharest

Faculty of Mathematics and Computer Science

Str. Academiei\ 14, 010014 Bucharest, Romania

E-mail: miculesc@yahoo.com, mihail\_alex@yahoo.com

\end{document}